\documentclass[12pt]{amsart}
\usepackage{amsmath}
\usepackage{amscd}
\usepackage{graphics}
\usepackage{latexsym}
\usepackage{amssymb}

\oddsidemargin .in
\theoremstyle{plain}
\newtheorem{Thm}{Theorem}

\newtheorem{Rm}{Remark}

\theoremstyle{definition}

\theoremstyle{remark}

\def\C{{\mathbb C}}

\def\R{{\mathbb R}}
\def\P{{\mathbb P}}

\def\Z{{\mathbb Z}}

\def\CM{\mathcal M}

\def\CT{\mathcal T}
\def\CF{\mathcal F}
\def\CA{\mathcal A}
\def\CB{\mathcal B}

\def\CZ{\mathcal Z}

\def\CZ{\mathcal Z}

\def\s2x{\hbox{$S^2 \times S^2$}}

    \def\sqr#1#2{{\vcenter{\hrule height.#2pt
            \hbox{\vrule width.#2pt height#1pt \kern#1pt
            \vrule width.#2pt}\hrule height.#2pt}}}
    \def\square{\mathchoice\sqr67\sqr67\sqr{2.1}6\sqr{1.5}6}

\def\qed{~\hfill$\square$}

\begin{document}

\title[]{ Real algebraic structures }
\author{Selman Akbulut }
\thanks{The author is partially supported by NSF grant DMS 9971440}
\address{Department  of Mathematics, Michigan State University,  MI, 48824}
\email{akbulut@math.msu.edu }
\date{\today}
\begin{abstract}
A brief survey of real algebraic structures on topological spaces  is given.
\end{abstract}
\maketitle

\setcounter{section}{-1}
\vspace{-0.4in}
\section{\bf Introduction}


The question of when a manifold $M$  is homeomorphic (or diffeomorphic) to a real algebraic set $V$ is an old one. If  we start with an imbedding $M\subset {\R}^n$ and insist on finding an algebraic subset $V$ of $\R^n$ which is isotopic to $M$ in $\R^n$,  the problem  encounters  additional difficulties coming from complexification. Hence it is natural to break the question into two parts: (1) Stable:  If $M$ homeomorphic (or diffeomorphic) to some real algebraic set.  (2)  Ambient: If $M$ isotopic to a real algebraic subset in $\R^n$. While the first problem has a complete solution, the second one has many interesting obstructions. Here we give a quick summary of some of the related results. Clearly this brief survey is by no means complete, and it is biased towards author's interest in the field. For the basics  reader can consult to the book \cite{ak3}.

\section {\bf Stable Results}

By \cite{n} and \cite{t} every closed smooth manifold is diffeomorphic to a nonsingular real algebraic set, and by \cite{ak5} every closed PL manifold is homemorphic to real algebraic set. Also if $M\subset V$ is a closed smooth submanifold of a nonsingular variety $V$, we can ask whether $M$ can be made algebraic in $V\times {\R}^k$ for some large k . The complete solution of this is also known, to explain this we first need some definitions: We call a homology cycle of $V$ {\it algebraic}  if it is represented by a real algebraic subset. For example,  it is known that all the Steifel-Whitney and Pontryagin classes of $V$ are represented by algebraic cycles \cite{ak1},\cite{ak8}.  Let  $H_{*}^{A}(V,\Z_{2})$ be the subgroup of $H_{*}(V,\Z_{2})$ generated by the real algebraic subsets.  We call a real algebraic set  {\it totally algebraic} if  $H_{*}(V,\Z_{2})= H_{*}^{A}(V,\Z_{2})$.  Cleary $\R^n$ is totally algebraic, and  the (unoriented) Grassmannians $G(k,n)$ of $k$ planes in $\R^n$ are  totally algebraic \cite{ak1} (because its homology is generated by the Schubert cycles which are algebraic subsets).

{\Thm (\cite{ak1}) Every closed smooth submanifold $f:M\hookrightarrow V$ of a real algebraic set $V$ is isotopic to a real algebraic subset in $V\times \R^k$ for some large $k$ if and only if $V$ is totally algebraic.
\[
\begin{array}{ccc}
  &   &  V\times{\R}^{k}  \\
  &\nearrow   &  \downarrow \\
M&\stackrel{f}{\hookrightarrow}   & V  
\end{array}
\]

} 
{\Rm It should be noted that not every closed smooth manifold $V$ is diffeomorphic to a nonsingular totally algebraic  set  \cite{bd1}, but surprisingly every closed smooth $V$ is homeomorphic to a totally algebraic set $Z$  \cite{ak4}. Evidently the singularities of $Z$ is related to the obstructions arising from}  \cite{bd1}.

\section{\bf Ambient Results}

{\Thm (\cite{ak6}) Every closed smooth submanifold $M\subset \R^n$ is $\epsilon$ -isotopic to nonsingular points of an algebraic subset $V\subset \R^n$. That is $M$ is isotoping to a topological component $V_0$  of a  real algebraic set 
$V$ which is nonsingular, and the other components $V-V_0$ are smaller dimensional.}

\vspace{.1in}

Here $\epsilon$ - isotopy means arbitrarily small isotopy. In the proof of the above theorem one can arrange so that the extra components $Z=V-V_0$ are arbitrarily far away from $V_{0}$. Also this theorem implies that any $M\subset {\R}^n$ is $\epsilon$ - isotopic to a nonsingular algebraic subset $Z\subset \R^{n+1}$. This is because if $V=f^{-1}(0)$ and $V-V_{0} =g^{-1}(0)$ for some polynomials $f(x)$ and $g(x)$, then we can take 
 $$Z= \{ (x,t) \;|\: f(x)=0 \; \; tg(x)=1 \}$$ 

There is also a more sophisticated version of this theorem for immersed submanifolds, to explain it we need some definitions: Recall that the Thom construction associates an imbedding of any closed smooth manifold
$f: M^{m}\hookrightarrow \R^{m+k}$ and element in the homotopy group of the Thom space 
$[f]\in\pi_{m+k}(MO_{k})$, which classifies imbeddings up to cobordism in  $\R^n \times [0,1] $. In the more general case of  $f$ is an immersion,   Thom construction gives an element in the homotopy group of the iterated loop space suspension  $[f]\in\pi_{m+k}(\Omega^{\infty}\Sigma^{\infty} MO_{k})$; this is a group which classifies immersions up to immersed cobordisms in  
$\R^n \times [0,1] $. Let $\pi_{m+k}^{alg} (\Omega^{\infty}\Sigma^{\infty} MO_{k})$ be the subgroup generated by the cobordism classes of `almost nonsingular' algebraic subsets. We call an algebraic subset of $ \R^n $  an {\it almost nonsingular algebraic subset} if it is an image of a smoothly immersed manifold $f: M  \looparrowright {\R}^n $,  and each sheet of the immersion is nonsingular. In particular if $f$ is an imbedding, then the almost nonsingular algebraic set corresponding to $f$ is a nonsingular algebraic subset. Now we can state:

{\Thm (\cite {ak7})An immersed closed smooth manifold $f: M \looparrowright {\R}^n$ is $\epsilon$- isotopic to an almost nonsingular algebraic subset of $\R^n$ if and only if 
$\;[f]\in \pi_{m+k}^{alg} (\Omega^{\infty}\Sigma^{\infty} MO_{k})$.} 

\vspace{.1in}
In particular, if an imbedding of a closed smooth manifold $f:M \hookrightarrow {\R}^{n}$ is cobordant through immersions to a closed smooth submanifold  $N\subset \R^{n-1}\times \{0\}$, then it is isotopic  to a  nonsingular algebraic subset. This is because, by the remarks following Theorem 2, $N$ is $\epsilon$-isotopic to a nonsingular algebraic subset of $\R^n$, then by Theorem 3  we can isotope $M$ to a nonsingular algebraic subset of $\R^n$. So the  relevant topological question is whether we can cobord $M\subset \R^n$ into a lower dimensional subspace of $\R^n$: 

\vspace{.05in}

(i) When the normal bundle of $M\subset {\R}^n$ splits a trivial line bundle?

\vspace{.05in}

(ii) When $ [f] $ lies in the image of the suspension map $\Sigma$?
$$\pi_{m+k-1}(MO_{k-1})\stackrel{\Sigma}{\to} \pi_{m+k}(MO_{k})$$

\noindent Answering the first question  would give conditions when $M$ immerses into $\R^{n-1}$, answering the second second question would help us decide if $M$ is cobordant to a submanifold of $\R^{n-1}$. Ideally one would hope to reduce the answers to the conditions on characteristic classes of the normal bundle of $M$, then (ii) would be sufficient conditions isotoping $M\subset \R^n$ to a nonsingular real algebraic subset.


\section {\bf Real algebraic characteristic numbers }

\vspace{.1in}


Surprisingly, we still don't know whether a closed smooth submanifold $M\subset {\R}^{n}$ is isotopic to a (singular or nonsingular) real algebraic subset (though we came close to answering this  in the affirmative in  Section 2). So one can try to find obstructions. Currently this can only be achieved either by relaxing the condition of smoothness of $M$ (this section), or  by strengthening the notion of non-singularity (Section 4). Underlying topological space of every algebraic set is a stratified space (a polyhedron in particular), so it is natural to generalize the above questions from smooth manifolds to stratified subspaces of $\R^n$. 

{\Thm (\cite{ck},}\cite{ak3},\cite{ak10}) {\it There is a stratified space $Z^{3}\subset \R^4$ which is homeomorphic to an algebraic set, but can not be isotopic to an algebraic subsets of $\R^4$.}

\vspace{.1in}
To explain this, we need to review the general program of topologically characterizing  real algebraic sets, in particular we need to recall some topologically defined structures (`topological resolution tower structures') on stratified spaces, which enable us to identify the obstructions of  making stratified spaces homeomorphic to algebraic sets. These structures on stratified sets give a topological model for algebraic sets.
 
\vspace{.1in}

In \cite{ak3} a topological characterization program for real algebraic sets is introduced, here is a brief summary: Just like an algebraic number such as $2+\sqrt {3}$ is determined by a pair of integers and a monomial $\{2,3, x^2=2\}$ (integers are glued by a monomial), an algebraic set $V$ is  determined by a collection of nonsingular algebraic sets and a collection of  ``compatible'' monomial maps between them $\CF= \{ V_{i}, p_{ij} \}$, and $V$ is obtained by gluing these nonsingular $V_{i}$'s by monomial $p_{ij}$'s. We denote this by $|\CF|$
\begin{equation}
 V= | \CF| = \cup {V_{i}}/p_{ij}(x)\sim x
 \end{equation}
(e.g. Figure 2 is obtained by gluing an $S^2$ and two copies of $S^1$ together by a fold map)
This result is obtained by resolving various strata of $V$. These objects $\CF= \{ V_{i}, p_{ij} \}$ are called {\it algebraic resolution towers}. By imitating this, we can define analogous objects in the topological category $\CF_{top}=\{ M_{i}, p_{ij} \}$, where $M_{ij}$ are smooth manifolds and $p_{ij}$ are certain topological version of monomial maps between them.  
As in (1),  by gluing the objects of $\CF_{top}$ by its maps we obtain certain stratified spaces $X=|\CF_{top}|$, called {\it topological resolution towers}. We have the following categories of sets:

\vspace{.05in}

\hspace{.15in}  $\CA\; =\;\{\CF\}$ \;\;\;\;\; \;\;: \; Algebraic resolution towers
   
 \hspace{.15in}  $\CT\;=\;\{\CF_{top}\}$ \;\;\;\;\;: \; Topological  resolution towers 
  
\hspace{.15in}  $ |\CA|\;=\;\{\;|\CF|\;\}$\; \;\;: \;  Realization of algebraic resolution towers 

\hspace{.15in}  $ |\CT|\;=\;\{\;|\CF_{top}|\;\}$ : \;\;Realization of topological resolution towers 
 
 \vspace{.05in}
We have the following  maps, where the vertical arrows are gluings (and they are onto by construction), the horizontal right-pointing arrows are the forgetful maps, and the bottom left-pointing arrow is the important {\it algebraization map}: it is a generalized version of Theorem 1 (it turns a collection of smooth manifolds and  compatible topological maps between them into nonsingular algebraic sets and compatible rational maps, such a way that gluing them gives an algebraic set).
\begin{equation}
\begin{array}{ccc}
       |\CA|   & \to & | \CT |\\
   \uparrow  & & \uparrow\\
     \CA   &\stackrel{ \longrightarrow }{\leftarrow} & \CT
\end{array}
\end{equation}
Therefore the stratified sets in  $|\CT |$ topologically characterizes the  real algebraic sets in $|\CA| $. Also  if  $\CA lg$  denotes the category of all real algebraic sets, by the above gluing process (1) we get a surjection $\gamma: \CA lg \to  |\CA|  $.  This appears to give a complete topological characterization of all real algebraic sets. Unfortunately this is not quite so.  For this we need that  the elements  $\CF=\{ V_{i}, p_{ij} \}$  in the image of $\gamma$  to be {\it submersive},  that is  we need the  maps  $p_{ij}$ to be ``submersive" on each strata. We may impose this property as the part of the definition of $\CT$ (this makes the diagram (2) commutative). This property is known to exits for  algebraic sets of dimension $ < 4$, and it would hold  in general provided that there is a certain map version of the ``resolution of singularities theorem'',  which  
we don't know if exists or whether it follows from a modification of the usual resolution of singularities theorem of Hironaka \cite{ak3}, \cite{k}.

\vspace{.05in}

There is an interesting subclass of stratified spaces in $|\CT|$,  called  $A$-{\it spaces}, which behave nicely on $PL$ manifolds. For example, their existence on $PL$ manifolds can be reduced to an algebraic topology problem, i.e. a bundle lifting problem (fortunately with zero obstruction). Furthermore $A$- spaces are submersive elements of $\CT$. This is why all $PL$ manifolds are homeomorphic to real algebraic sets. 

{\Thm (\cite{ak2} \cite{ak5}, \cite{at}) \it Every closed $PL$-manifold is homeomorphic to an algebraic set.}

\vspace{.1in}

Therefore any topological obstruction for a stratified space $X$ to lie in $|\CT|$ is an obstruction $X$ to be isomorphic (as stratified space) to an algebraic set.  It is already known by Sullivan, that every real algebraic set must be an {\it Euler space} (a stratified set such that the link of every point has even Euler characteristic).  It turns out that in dimensions $\leq 2$   this is also the sufficient condition for a stratified set to be homeomorphic to an algebraic set. This is proven by showing that every $2$-dimensional Euler spaces lie in $|\CT |$ ([AK1-3], \cite{bd2}). By studying the topology of  $\CF_{top}=\{ M_{i}, p_{ij} \}$ carefully, one can start defining inductively a sequence of characteristic numbers on $n$-dimensional Euler spaces $X_{n}$  (here $\Z_{2}=\Z/2\Z$) 
$$\beta=\beta(n):X_{n}^{(0)}\to  \Z_{2}^{ d(n)}$$  
whose vanishing is necessary and sufficient (or just necessary) for $X_{n}$ to lie in $|\CT|$, where $X_{n}^{(0)}$ is the $0$-skeleton of $X_{n}$. In our context this means that the $n-1$ dimensional links  of the verticies of $X$, which are already in $\CT$, should bound in $\CT$. Roughly $\beta_{n}$ are  the the cobordism characteristic numbers of $n-1$ dimensional elements of $\CT$. 

\vspace{.1in}

In \cite{ak3} this program was carried out in the first nontrivial case  of $n=3$. It turns out $d(3)=4$.  In fact it was shown that $X_{3}$ is homeomorphic to an algebraic set if and only if its characteristic numbers $\beta =0$. It goes as follows: To every $1$-dimensional Euler space $X_{1}$ we associate numbers $\alpha_{j}= \alpha_{j} (X_{1})$ $j=0, .. , 7$,  which are number of vertices of $X_{1}$ whose links has $j \; (mod \;8) $ number of points. Then to $X_{1}$  we associate  the following (well defined after subdivision) $4$-tuple number  (e.g. Figure 1).

  \begin{figure}[ht]  \begin{center}  
\includegraphics{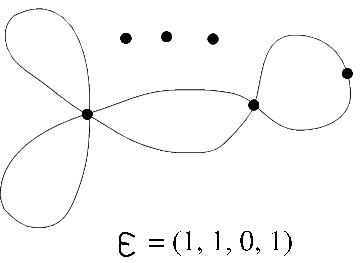}   \caption{}    \end{center}
   \end{figure} 

\vspace{-.1in}
$$\epsilon(X_{1})=(\alpha_{0}, \alpha_{6}, (\alpha_{0} + \alpha_{4})/2, (\alpha_{2}+\alpha_{6})/2) \in \Z_{2}^{4}$$

\noindent Hence to the verticies of any  $2$-dimensional Euler space $X_{2}$ we can associated $4$-tuple numbers (because their links are $1$-dimensional Euler spaces). For the example of Figure 2 we calculated these numbers (and drew  the links of its verticies). 
 \begin{figure}[ht]
 \begin{center}  
\includegraphics{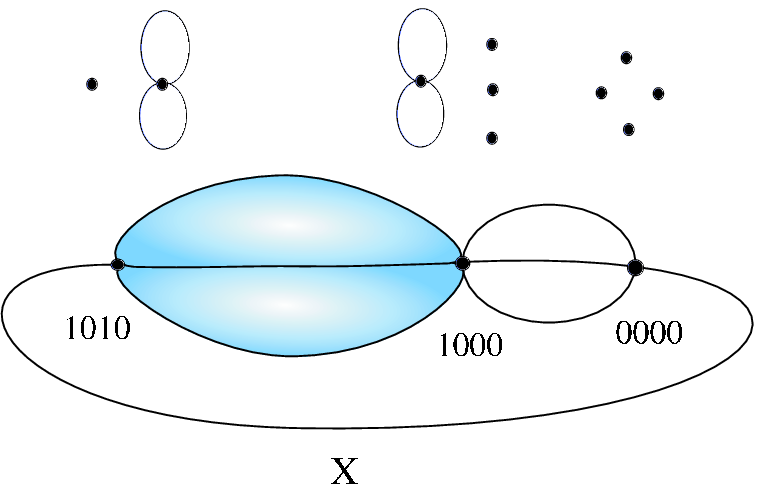}  \caption{}   
 \end{center}
   \end{figure}

Now for  $\epsilon\in \Z_{2}^{4}$, we can associate numbers to the verticies of any $3$-dimensional Euler space: $\beta_{\epsilon} :X^{(0)}_{3}\to \Z_{2}$  by  
 $p \mapsto ${\it number of vertices (mod 2) in $L:=link(p)$,  with $\epsilon (L) =\epsilon$}. Then the definition of $\beta :X^{(0)}_{3}\to \Z_{2}^4$  is given by the following expression
$$ \beta= (\beta_{0100}+  \beta_{0101}, \;  \beta_{1000}+  \beta_{1001}, \; \beta_{1100}+  \beta_{1101}, \;  \beta_{1110}+  \beta_{1111})$$

For example, if $X_{3}$ is the suspension of the $2$-dimensional Euler space of Figure 2, it is not homeomorphic to an algebraic set; while  the suspension $Z$ of the  Euler space in Figure 3 is homeomorphic to  a real algebraic set. This $Z$ has the property mentioned in Theorem 4.
\begin{figure}[ht]  \begin{center}  
\includegraphics{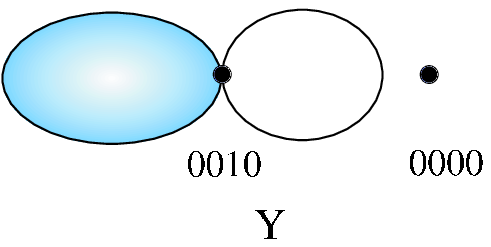}   \caption{}    \end{center}
   \end{figure}
   
    More recently the ``necessary" part of this program was generalized to dimension $4$  by \cite{p} and \cite{mp2} (where \cite{p} generalizes this approach, whereas \cite{mp2} uses an alternative approach),  and also  \cite{p} finds sufficient invariants (not all in $\Z_{2}$ but in a larger group). One difficulty with this is $d(n)$ grows too fast to calculate, for example \cite{mp2} shows $d(4)\geq 2^{43}-43$.

 \vspace{.1in}

 Note that in examples above the Euler characteristic $mod$ $4$ of the link of the $1$-dimensional strata is generically constant. This is no coincidence, in fact it was observed in \cite{ck} that general  real algebraic varieties  $V\subset X$ satisfy the property (*):  The Euler characteristic $\chi (lk_{x}(V,X)) $ of the link of $V$ in $X$ ($ x\in V$)  is generically constant  $mod \;4$ (also see \cite{ak3},  \cite{mp1}). In \cite{ck}  this property is used  to show that the set $Z\subset \R^4$ above satisfies the second claim of Theorem 4 as follows:  Let $I\subset Z$ be the arc, which is the suspension of the isolated point of $Y$.  By (*)  in any algebraic model of $Z$, the Zariski closure of $I$ is a $1$-dimensional algebraic subset $V\subset Z$ such that  $ I'=V-I $ generically lies in the $3$-manifold part of $Z$. Now if  $Z\subset \R^4$  were an algebraic hyper-surface, we can choose a square free polynomial equation $f(x)=0$ of $Z.$ Say,  $f$ takes signs $\pm$  on the inside and outside regions $B_{\pm}\subset \R^4$ separated by $S^3\subset Z$. Define
   $$\widetilde{\R^4}:=\{(x,t)\in \R^4\times \R \;| \; t^2=f(x)\} \supset Z\supset V $$ 
 Then $\chi(lk_{x}(V,\widetilde{\R^4}))-\chi(lk_{x}(V,\R^4))=\chi (lk_{x}(V,B_{+}))-\chi (lk_{x}((V,B_{-}))$ mod $4$ is not generically constant (it is $2$ or $0$ when  $x$ is in the interior of $I$ or $I'$), this violates (*).

\section {\bf Transcendental manifolds}

If we ask whether a smooth submanifold $M\subset \R^n$ is isotopic to a nonsingular real algebraic set $V$ in a strong sense, then we can find genuine obstructions to doing this, even when $M$ is already nonsingular algebraic set in $\R^n$. Here nonsingular in the ``strong sense'' means $V$ is the real part of a nonsingular complex variety in $\C\P^n$, where $\R^n\subset\R\P^n$  is identified with one of the  standard  charts.

{\Thm (\cite{ak9}) There are smooth submanifolds of $M\subset \R\P^n$ which are isotopic to nonsingular real algebraic subsets, but can not be isotopic to the real parts of nonsingular complex algebraic subsets of $\C\P^{n}$}.
\vspace{.1in}

We will break the proof into small elementary steps which are mostly special cases of more general results  (e.g. \cite{ak3}), some of which are already mentioned in Section 1. When possible we will outline the proofs from scratch for  the benefit of non-specialist. 

\vspace{.1in}

$\bullet$ {\it {\bf Step 1}: Grassmannians of $k$-planes in $\R^n$ is a nonsingular real algebraic set: 
$$ G_{k}({\R}^n)=\{A\in {\R}^{n(n+1)/2}\;|\; A^{2}=A, trace(A)=k\;\}$$
where ${\R}^{n(n+1)/2}$ denotes the $n\times n$ symmetric matrices.  
If $V^{n-k}\subset {\R^n}$ is a compact 
nonsingular real algebraic set, the normal (or tangent) 
Gauss map $\gamma :V\to  G_{k}({\R}^n)$ 
is an entire rational map (non-zero denominator). Hence in particular the duals of  Steifel-Whiney classes of $V$ are represented by real algebraic subsets (since Steifel-Whitney classes of the universal bundle is represented by Shubert cycles  in $G_{k}({\R}^n) ).$}

\proof We can identify $G_{k}({\R}^n)$ with the algebraic  set defined on the right via the map $L\mapsto \mbox{projection matrix to\;} L$. Since nonsingular $V$ can be covered by finitely many Zariski-open sets $V=\cup_{j=1}^{N} V_{j}$,  each of which is described in $\R^n$ as the zeros of $k$ polynomials  $f_{i}(x)=0$,  $i=1,.., k$,  whose Jacobian $ A(x)=(\partial{f_{i}}/\partial{x_{j}})$ has rank $k$. Then by standard linear algebra $\gamma(x)= A(x)(A^{t}(x)A(x))^{-1}A^{t}(x)$. In particular this allows us to describe the restriction of the Gauss map  
$\gamma|_{V_{j}}=P_{j}(x)/q_{j}(x)$ as a matrix whose entries are entire rational maps. Then $\gamma =\sum _{j=1}^{N} P_{j}/q_{j}$ is the answer. \qed 

\vspace{.15in}

$\bullet$ {\it {\bf Step 2}:  If $V\subset \R\P^{N}$ is a compact nonsingular real algebraic set with a nonsingular complexification $j: V\hookrightarrow V_{\C}\subset \C\P^{N}$,  and let $L\subset V$ be an algebraic subset with complexification $L_{\C}\subset V_{\C}$. Then the restriction  of the Poincare dual of the fundamental class of $L_{\C}$ is the cup square of the Poincare dual of the fundamental class of $L$, i.e.  as $\Z_{2}$ classes:}$$j^{*}PD[L_{\C}]=PD[L]^{2}$$ 
\proof Let $g_{\C} :\tilde{L}_{\C}\to V_{\C}$ be the resolution of singularities map
 followed by the inclusion $(\tilde{L}_{\C}, \tilde{L}) \to (L_{\C}, L)\hookrightarrow (V_{\C}, V) $.  Isotop $g_{\C}$  to a map 
$g' :\tilde{L}_{\C}\to V_{\C}$ which is transverse to $V$, then  $j^{*}[L_{\C}]$ is the Poincare dual of the   intersection $g' (L_{\C})\cap V$. Call  $g_{\C}|_{\tilde{L}}=g$. Near $V$ the map $g_{\C}$ is modeled by the map $g\times g: (\tilde{L}\times \tilde{L}, \Delta) \to (V\times V, \Delta)$ where $\Delta$ are the diagonals, i.e. $(x,y)\mapsto (g(x),g(y))$. Therefore $j^{*}[L_{\C}]$ is represented by the Poincare dual of the self intersection of the class $[L]$ \qed

\vspace{.1in}

$\bullet$ {\it  {\bf Step 3}: If  $V\subset {\R^n}$ is an algebraic set given by polynomials with highest degree 
terms are $|x|^{2d}$, then $V$ is projectively closed. That is  if $\lambda: \R^n\hookrightarrow \R\P^n$ is 
the imbedding  $(x_1,..,x_n)\mapsto [1,x_1,..,x_n]$, then $\lambda(V)$ is a projective 
algebraic set in $\R\P^n$.}

\proof If $f(x)=0$ is a defining polynomial equation of $V$ (by taking the sum of the squares of the defining equations, every real algebraic set in $\R^n$ can be described by a single polynomial equation). Let $f_{d}(x)$ be the highest degree term of $f(x)$. Then clearly the equations of $\lambda(V)$ is $0= x_{0}^{d}f(x/x_{0})=f_{d}(x) +x_{0 }q(x)$, for some polynomial $q(x)$,  hence the zeros $\lambda(V)$  in $\R\P^n$ coincides with the zeros of $V$. \qed

\vspace{.15in}

$\bullet$ {\it  {\bf Step 4}: If $M^{m}\subset Y^{m+1}\hookrightarrow \R^n$ are imbeddings of closed smooth manifolds, such 
that $M$ is separetes $Y$, then $M$ is $\epsilon$ - isotopic to a projectively closed nonsingular 
algebraic subset of  $\R^n$.}

\proof Let $Z^{n-1}\subset \R^n$ a codimension one smooth submanifold of which intersects $Y^{n+1}$ transversally at  $M^{m}$, e.g. we can take $Z$ to be the boundary of narrow sausage  as shown in the Figure 4.  
By Theorem 2 we can approximate $Y$ by a component $Y'_{0}$  of a nonsingular  of an algebraic set $Y'$. We claim that $Z^{n-1}$ can be approximated by a  projectively closed nonsingular algebraic subset $Z' \subset \R^n$.

  \begin{figure}[ht]  \begin{center}  
\includegraphics{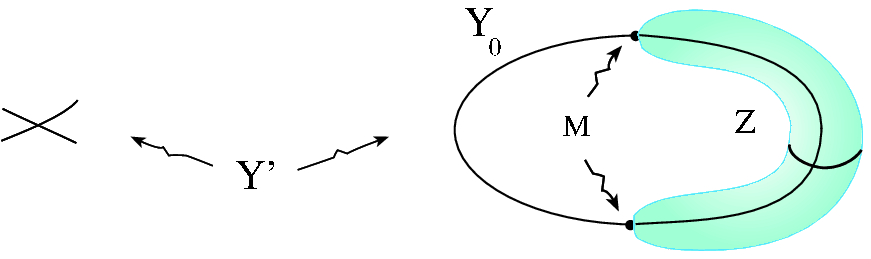}   \caption{}    \end{center}
   \end{figure}

This is because we can write  $Z=f^{-1}(0)$, where $f: \R^n\to \R$ is some  proper smooth function with $0$ as a regular value (since $Z$ is codimension one), and then approximate $f$ by a polynomial of the form $F=g +\epsilon |x|^{2k}$ where $k>>1$ and take $Z'=F^{-1}(0)$. Therefore by Step 3 the set $Y'\cap Z'$  is a projectively closed nonsingular algebraic set isotopic to M. \qed

\vspace{.15in}

\proof (of Theorem 5) We can construct examples from imbeddings $\R\P^{m}\subset \R^{2m-s}$. For a nonorientable example we need  $m$ even and $s\geq 3$, and for an orientable example we need $m=4k+1$ and $s\geq 5$. For example \cite{nu} gives $\R\P^{10}\subset \R\P^{17}$, and \cite{th} gives $\R\P^{13}\subset \R^{20}$. Let us do the nonorientable case: Take the obvious imbeddings
$$ M^{11}=\R\P^{10}\times S^{1} \subset \R^{17}\times \R^{3}=\R^{20} \subset  \R\P^{20}$$
Let $Y^{12}=\R\P^{10}\times S^{2}$ in  $\R\P^{20}$. Here $S^{1}\subset S^{2}$ are the standard spheres. By Step 4 we can isotope $M^{11}$ to a nonsingular projective algebraic subset $V^{11}\subset \R\P^{20}$. We claim that $V^{11}$ can not be the real part of a nonsingular complex algebraic subset $V_{\C}$ of  $\C\P^{20}$ (defined over $\R$).  Suppose such a $V_{\C}$ exists.
 By the Lefschetz hyperplane theorem (e.g. \cite {h}),
 for $i\leq 2(11)-20=2$  the  restriction induces an 
isomorphism
$$ H^{i}(\C\P^{20};\Z)\stackrel{\cong}{\to} H^{i}(V_\C;\Z)$$
By Step 1 the Poincare dual the Steifel-Whitney class $w_{1}(V)=\alpha \times 1$ is represented by an algebraic subset $L\subset V$; here $\alpha \in H^{1}(\R\P^{10};\Z_{2}) $ is the generator.   But since $V$
lies in a chart $\R^{20}$ of $\R\P^{20}$, the right vertical (restriction) map in the following commuting diagram is the zero map. This is a contradiction to Step 2.
    $$\begin{array}{ccccc}
PD[L_{\C}] \in  \hspace{-.1in}&H^{2}(V_\C; \Z _{2}) & \stackrel{j^{*}}{\longrightarrow } 
    & H^{2}(V; \Z _{2})& \hspace{-.1in} \ni PD[L]^{2} =\alpha^{2}\times 1 \neq 0\\
&{ \small \cong }\uparrow &                   &  \uparrow \mbox{\small zero} &\\
 & H^{2}(\C\P^{20}; \Z _{2}) 
  & \stackrel{j^{*}}{\longrightarrow } & H^{2}(\R\P ^{20}; \Z _{2}) &
\end{array}$$

\vspace{.1in}

For an orientable example we start with  $\R\P^{13}\subset \R^{20}$ and use $w_{2}(\R\P^{13})\neq 0.$
\qed                                                                        

{\Rm Though the likelihood is slim, If one can demonstrate an imbedding of a smooth $M\subset \R^{n}$ such that one of its Steifel-Whitney or Pontryagin  classes  can only be represented by a bad singular space $X\subset M$ with the property  of Theorem 4,  then by Step 1  $M$ can not be isotopic to a nonsingular algebraic subset of $\R^n$.}

{\Rm Recall in  gauge theory, the space of  connections on an $SU(2)$-bundle $P\to M^4$ modulo the gauge group of $P$,  is identified by a component of the mapping space $\CB(P)= Map(M,G_{3}(\R^n))$, where $n>>1$. Prescribing  a metric $g$ on $M$ allows us to consider the space self-dual connections $\CM_{g} \subset \CB(P)$ which is finite dimensional. Then evaluating various cohomology classes of $\CB(P)$  on $\CM_{g}$ gives Donaldson invariants. So just as using metrics on $M^4$ as auxiliary objects to define natural finite dimensional subspaces $\CM_{g}$, one can try to use real algebraic structures on a smooth submanifold $M^{n-k}\subset \R^{n}$ as auxiliary objects to define finite dimensional subsets $\CZ_{d} \subset Map(M,G_{k}(\R^n))$, such as $\CZ_{d}=\{$entire rational maps of degree $\leq d \}$,  and imitate Donaldson invariants (here compactness is the main problem to be faced).}

 \end{document}